\documentclass[number,citesort,seceqn,dvips]{arxbj}
\usepackage{upgreek}

\aid{0}
\volume{17}
\issue{3}
\pubyear{2011}
\firstpage{880}
\lastpage{894}
\doi{10.3150/10-BEJ297}

\makeatletter

\newproclaim{remark}{Remark}
\newproclaim{remarks}{Remarks}

\def\eps{\varepsilon}

\newcommand{\eqref}[1]{(\ref{#1})}

\makeatother

\begin{document}
\begin{frontmatter}

\title{Renorming divergent perpetuities}
\runtitle{Renorming divergent perpetuities}

\begin{aug}
\author[a]{\fnms{Pawe{\l}} \snm{Hitczenko}\corref{}\thanksref{a}\ead[label=e1]{phitczenko@math.drexel.edu}%
\ead[url,label=u1]{http://www.math.drexel.edu/\textasciitilde phitczen}}
\and
\author[b]{\fnms{Jacek} \snm{Weso{\l}owski}\thanksref{b}\ead[label=e2]{wesolo@mini.pw.edu.pl}%
\ead[url,label=u2]{http://www.mini.pw.edu.pl/\textasciitilde wesolo}}
\runauthor{P. Hitczenko and J. Weso{\l}owski}
\address[a]{Departments of Mathematics and Computer Science,
Drexel University,
Philadelphia, PA 19104,
USA.
\printead{e1},
\printead{u1}}
\address[b]{Wydzia{\l} Matematyki i Nauk Informacyjnych,
Politechnika Warszawska,
Plac Politechniki~1,
00-661 Warszawa, Poland.
\printead{e2},
\printead{u2}}
\end{aug}

\received{\smonth{8} \syear{2009}}
\revised{\smonth{3} \syear{2010}}

%
\begin{abstract}
We consider a sequence of random variables $(R_n)$ defined by the
recurrence $R_n=Q_n+M_nR_{n-1}$, $n\ge1$, where
$R_0$ is arbitrary and $(Q_n,M_n)$, $n\ge1$, are i.i.d. copies of a
two-dimensional random vector
$(Q,M)$, and $(Q_n,M_n)$ is independent of $R_{n-1}$.
It is well known that if $E{\ln}|M| <0$ and $E{\ln^+}|Q|<\infty$, then
the sequence $(R_n)$ converges in distribution to a random variable $R$
given by $R\stackrel d=\sum_{k=1}^\infty Q_k\prod_{j=1}^{k-1}M_j$, and
usually referred to as perpetuity. In this paper we consider a~%
si\-tuation in which the sequence $(R_n)$ itself does not converge. We
assume that $E{\ln}|M|$ exists but that it is non-negative and we ask if
in this situation the sequence $(R_n)$, after suitable normalization,
converges in distribution to a non-degenerate limit.
\end{abstract}

%

\end{frontmatter}

\section{Introduction}
We consider the following iterative
scheme
%
\begin{equation}\label {perp-nq}
R_n=Q_n+M_nR_{n-1},\qquad  n\ge1,
\end{equation}
where
$R_0$ is arbitrary and $(Q_n,M_n)$, $n\ge1$, are i.i.d. copies of a
two-dimensional random vector
$(Q,M)$, and $(Q_n,M_n)$ is independent of $R_{n-1}$.
Writing out
the above recurrence we see that
%
\begin{equation}\label {perp-exp}
R_n=Q_n+M_nQ_{n-1}+M_nM_{n-1}Q_{n-2}+\cdots+M_n\cdots M_2Q_1+M_n\cdots
M_1R_0.
\end{equation}
Note that although $(R_n)$ is not a sequence of partial sums, after
renumbering $(Q_n,M_n)$ in the opposite direction, we can write
%
\begin{equation}
R_n\stackrel d=\prod_{j=1}^nM_jR_0+\sum_{k=1}^nQ_k\prod_{j=1}^{k-1}M_j,
\label {perp-it}
\end{equation}
where \lq\lq$\stackrel d=$\rq\rq\ denotes the equality in
distribution and we adopt the convention that the product or sum over
the empty range is 1 or 0, respectively. Much of the impetus for
studying such equations stems from numerous applications of schemes
like (\ref{perp-nq}) in mathematics and other disciplines of science.
We refrain here from giving a long list of fields in which equation
(\ref{perp-nq}) appeared, referring instead to some of the references
we give (most notably to \cite{eg,vervaat} for the status up to the
early nineties of the past century) for more detailed information.
Examples of more recent applications closer to statistics are given in
\cite{br}, an application related to neuronal modeling may be found in
\cite{hm} and, for an application in the analysis of algorithms, see
\cite{g} and references therein.

Most of the up-to-date research focused on the situation when the
sequence $(R_n)$ converges in distribution and on analyzing properties
of its limit. It has been known for a long time (see \cite{vervaat},
Lemma 1.7, or \cite{kesten}) that if
\[
E{\ln}|M| <0 \quad \mbox{and}\quad  E{\ln^+}|Q|<\infty,
\]
then the sequence of partial sums in (\ref{perp-it}) converges almost
surely and the product in the first term asymptotically vanishes. Thus,
$(R_n)$ converges in distribution to a random variable $R$ given by
\[
R\stackrel d=\sum_{k=1}^\infty Q_k\prod_{j=1}^{k-1}M_j.
\]
This $R$, referred to as perpetuity, is of central interest. In
particular, its tail behavior has been thoroughly investigated (see,
e.g., \cite{goldie,gg,grey,hw,kesten,vervaat} and references therein
for more information).

In this paper we focus on a situation in which the sequence $(R_n)$
does not converge. We will refer in such situations to the whole
sequence $(R_n)$ as a divergent perpetuity. We will assume throughout that
$E{\ln}|M|$ exists but that
%
\begin{equation}\label{mu}
\mu:=E{\ln}|M|\ge0.
\end{equation}
Under the above assumption we ask if $(R_n)$ can be renormed to
converge in distribution to a non-degenerate limit. As is the case when
$(R_n)$ converges, the roles of $R_0$ and $Q$ seem to be of much less
importance than the role of $M$. In fact, since $R_0$ plays no
significant role from now on we assume that $R_0=0$. Furthermore,
assuming that $M$ is non-random would lead to $(R_n)$ being a sequence
of sums of independent random variables. Since this situation has been
extensively studied, we will exclude it from our considerations by
assuming from now on that $M$ is non-constant. On the other hand, we
freely impose moment conditions on $Q$ when necessary, or sometimes
even assume that it is of a special form. For example, if $Q=1$, then
(\ref{perp-it}) is just a partial sum of successive partial products
of i.i.d. copies of $M$. While there is very little known technical
connection, it is perhaps worth mentioning that an analogous problem of
investigating the asymptotic properties of the consecutive products of
partial sums
\[
\prod_{k=1}^n\sum_{j=1}^k M_j \quad \mbox{vs.} \quad \sum_{k=1}^n\prod
_{j=1}^k M_j,
\]
has been intensely investigated in the past several years; see, for
example, \cite{rw,q,hz}.

As $x\to\ln x$ is concave, assumption (\ref{mu}) implies that
necessarily $E|M|\ge1$. (We do not assume, however, that $E|M|$ is
finite.) Thus there are at least four situations to consider, namely,
\begin{longlist}[(iii)]
\item[(i)] $E{\ln}|M|>0$, $E|M|>1$, $|M|$ non-random,
\item[(ii)] $E{\ln}|M|>0$, $E|M|>1$, $|M|$ random,
\item[(iii)] $E{\ln}|M|=0$ and $E|M|>1$,
\item[(iv)] $E{\ln}|M|=0$ and $E|M|=1$.
\end{longlist}
These cases are discussed in detail below and we will show that in each
of these situations $(R_n)$ can be renormed (each time differently) so
that it converges in distribution to a non-degenerate limit. As was
brought to our attention by A.~Iksanov, some particular cases of (ii)
and (iii) were studied by Rachev and Samorodnitsky \cite{rs} who, for
non-negative $Q$ and $M$ and under the assumption that $\log M_n$
belongs to the $\alpha$-stable domain of attraction, obtained the
log-stable limit law for suitably normalized $(R_n)$ -- see also
comments on this connection in Sections~\ref{(ii)} and~\ref{(iii)}.
Furthermore, weak convergence in the situation complementary to \eqref
{mu} is considered in \cite{pakes} where, for non-negative $Q$ and
$M$, it is assumed that $-\infty<\ln M<0$, but the perpetuity is
divergent because $E\ln^+Q=\infty$.

\section{The case $E{\ln}|M|>0$, $E|M|>1$ and non-random $|M|$}
In this section we consider the following situation: for a fixed $\rho>1$ let
\[
M\stackrel d=\rho\eps,\qquad  \mbox{where } P(\eps=1)=p,\ P(\eps=-1)=1-p=q.
\]

\begin{thm}\label {prop:sym} Let $M$ be as above. Assume that
$Q\stackrel
d=\eps$ and that it is independent of $M$.
\begin{longlist}[(ii)]
\item[(i)] Symmetric case: Let $p=1/2$. Then
\[
\frac{R_n}{\rho^{n-1}}\stackrel{d}{\longrightarrow}\sum
_{k=1}^{\infty}\lambda ^{k-1}\eps_k ,
\]
where $\lambda =\rho^{-1}<1$ and $(\eps_k)$ is the sequence of i.i.d.
copies of $\eps$.
\item[(ii)] Asymmetric case: Suppose $p\ne1/2$. Then we have
%
\begin{equation}\label {eq:as}\frac{R_n}{\rho^{n-1}}\stackrel
d\longrightarrow
rX,
\end{equation}
where
$r$ is a symmetric Bernoulli random variable (i.e., $P(r=1)=P(r=-1)=1/2$),
$X\stackrel
d=\sum_{k=0}^\infty\lambda ^k\prod_{j=1}^k\eps_j$
and $r$ and $X$ are independent.
\end{longlist}
\end{thm}

\begin{remark*}
The sums of the form
$\sum_{k=1}^{\infty}\lambda ^{k-1}\eps_k $ considered in part (i) are
well-known objects, usually referred to as \lq\lq symmetric
Bernoulli convolutions\rq\rq. Their properties have been extensively
studied since
mid-1930s. In particular, it is known that the limiting distribution is
uniform on
$[-2,2]$ if $\lambda =1/2$, singular if $0<\lambda <1/2$ and absolutely
continuous for almost all (but not all)
$\lambda \in(1/2,1)$. A good description of the current state of
knowledge can be found, for example, in~\cite{pss}.
\end{remark*}

\begin{pf*}{Proof of Theorem~\ref{prop:sym}}
We have
%
\begin{equation}\label {eq:rad} R_n=\eps_n+\rho\eps_n\eps
_{n-1}+\cdots+\rho
^{n-1}\prod_{j=1}^n\eps_j\stackrel d=
\sum_{k=1}^{n}\rho^{k-1}\prod_{j=1}^k\eps_{j}.
\end{equation}
If $p=1/2$ then
the sequences
%
\begin{equation}\label {eq:dist_rad} (\eps_1,\eps_1\eps_2,\eps
_1\eps_2\eps
_3,\dots) \quad \mbox{and}\quad  (\eps_1,\eps_2,\eps_3,\dots)
\end{equation}
are identically distributed. Thus, if we normalize $R_n$ by
$\rho^{n-1}$, we get
\[
\frac{R_n}{\rho^{n-1}}\stackrel d=\sum_{k=1}^{n}\rho^{k-n}\eps_k
\stackrel d=\sum_{k=1}^{n}\lambda ^{k-1}\eps_k .
\]
Since $0<\lambda <1$ the series of partial sums converges almost surely
and thus $R_n/\rho^{n-1}$ converges in distribution to the given limit.

If $p\ne1/2,$ then
the distributional equality in (\ref{eq:dist_rad}) is no longer valid. However,
we can write the right-hand side of \eqref{eq:rad} as
\[
R_n\stackrel d=
\sum_{k=1}^{n}\rho^{k-1}\prod_{j=1}^k\eps_{j}=
\rho^{n-1}\prod_{j=1}^n\eps_j \Biggl(\sum_{k=1}^{n}\frac{\lambda
^{k-1}}{\prod_{j=1}^{k-1}\eps_j} \Biggr).
\]
Since $\eps_k=1/\eps_k$ we get
%
\begin{equation}\label {eq:tn}\frac{R_n}{\rho^{n-1}}\stackrel
d=T_n\prod
_{j=1}^n\eps_j ,
\end{equation}
where
\[
T_\ell=\sum_{k=1}^{\ell}\lambda ^{k-1}\prod_{j=1}^{k-1}\eps_j .
\]
Fix arbitrary $m$. Then, for any $n>m$, write the term on the
right-hand side of (\ref{eq:tn}) as
\[
T_m \prod_{j=1}^n\eps_j+\prod_{j=1}^n\eps_j \Biggl(\sum
_{k=m+1}^{n}\lambda
^{k-1}\prod_{j=1}^{k-1}\eps_j \Biggr).
\]
The second summand is bounded in absolute value by $\lambda
^m/(1-\lambda )$
and thus it can be made arbitrarily small by choosing $m$
sufficiently large. Consequently, we consider only the first part
of the above expression. Note that $T_m$ depends on
$(\eps_1,\ldots,\eps_m)$ only. Therefore upon conditioning on $(\eps
_1,\ldots,\eps_m)$ we obtain
\[
E (\mathrm{e}^{\mathrm{i}tT_m\prod_{j=1}^n \eps_j} )=E \Biggl[
\phi_{n-m} \Biggl(tT_m\prod_{j=1}^m \eps_j \Biggr) \Biggr],
\]
where $\phi_k$ is the characteristic function of the product
$\prod_{j=1}^k \eps_j$. Since $\prod_{j=1}^n\eps_j=1$ if and only
if the number of $j$'s such that $\eps_j=-1$ is even, we have
\begin{eqnarray} \label{eq:symB}
P\Biggl(\prod_{j=1}^n\eps_j=1\Biggr)-P\Biggl(\prod
_{j=1}^n\eps_j=-1\Biggr)
&=&\mathop{\sum_{k=0}}_ {k\mbox{-}\mathrm{even}}^n\pmatrix{n\cr
k}q^kp^{n-k}-\mathop{\sum_{k=0}}_{k\mbox{-}\mathrm{odd}}^n\pmatrix{n\cr k}q^kp^{n-k}\nonumber\\[-8pt]\\[-8pt]
&=&(p-q)^n,\nonumber
\end{eqnarray}
which vanishes as $n\to\infty$. Thus the product $\prod_{j=1}^n\eps
_j$ converges in distribution
to $r$ (recall that its characteristic function is $\cos(t)$).
Therefore for arbitrary fixed $m$
\[
E (\mathrm{e}^{\mathrm{i}tT_m\prod_{j=1}^n\eps_j} )\stackrel{n\to\infty
}{\longrightarrow}
E \cos\Biggl(tT_m\prod_{j=1}^m \eps_j \Biggr)=E \cos(tT_m) .
\]

As $m\to\infty$ the sum $T_m$ converges almost surely to $X$. By
the Lebesgue dominated convergence it follows therefore that
\[
E \cos(tT_m)\stackrel{m\to\infty}{\longrightarrow}E \cos(tX) ,
\]
the latter being the characteristic function of $rX$, with $r$ and $X$
independent. This proves~(\ref{eq:as}).\hfill\mbox{} \mbox{\quad}
\end{pf*}

\section{The case $E{\ln}|M|>0$, $E|M|>1$ and random $|M|$}\label{(ii)}
In this section we assume that $\mu=E{\ln}|M|>0$ and that $|M|$ is
random. This forces $E|M|>1$, but it may be infinite. Set
$v^2:=\operatorname{var}({\ln}|M|)$. Then we have
\begin{thm} \label {prop:rand|M|} Let $(R_n)$ be given by \eqref
{perp-it} with the pair $(Q,M)$ satisfying $\mu>0$, \mbox{$0<v^2<\infty$},
and $E{\ln^+}|Q|<\infty$.
\begin{enumerate}[(ii)]
\item[(i)] As $n\to\infty$,
\[
\frac{|R_n|^{1/(v\sqrt{n})}}{\exp(\mu\sqrt
n/{v})}\stackrel{d}{\longrightarrow} \mathrm{e}^{{\mathcal N}},
\]
where ${\mathcal N}$ is a standard normal random variable.
\item[(ii)] Assume $P(M>0)>0$ and $P(M<0)>0$ and define for real $t$
and $x$, $x^t:=\operatorname{sgn}(x)|x|^t$.
Then we have
\[
\frac{R_n^{1/{(v\sqrt{n})}}}{\exp(\mu\sqrt
n/{v})}\stackrel{d}{\longrightarrow} r \mathrm{e}^{\mathcal N},
\]
where $r$ is a symmetric Bernoulli random variable independent of
${\mathcal N}$.
\end{enumerate}
\end{thm}

\begin{remark*}
The first part of this theorem overlaps with
Theorem~2.1\textup{(a)} in \cite{rs} where, for non-negative $Q$ and $M$ and
under the assumption that $\ln M$ is in the domain of attraction of an
$\alpha$-stable law, $1<\alpha\le2$, the authors obtained \textup{(i)} with
$\mathcal N$ replaced by an $\alpha$-stable random variable and suitable
normalization of $(R_n)$ on the left-hand side (the normalizing
constants are the constants implied by the definition of a domain of
attraction).
\end{remark*}

\begin{pf*}{Proof of Theorem~\ref{prop:rand|M|}}
Consider $R_n$ given in \eqref{perp-exp} and factor the product of
$M_j$'s to write it as
%
\begin{equation}\label {eq:rand|M|}R_n
= \Biggl(\prod_{j=1}^{n}M_j \Biggr)
\sum_{k=1}^nQ_k\prod_{j=1}^{k}\frac1{M_j}.
\end{equation}
Consider the first factor. By the classical CLT
\[
\frac{|\prod_{j=1}^nM_j|^{1/{(v\sqrt{n})}}}{\exp(
(\mu /v)\sqrt n)}=
\biggl(\frac{\prod_{j=1}^n|M_j|}{\mathrm{e}^{\mu n}} \biggr)^{1/{(v\sqrt{n})}}
= \exp\biggl\{ \frac{\sum_{j=1}^n{\ln}|M_j| -n\mu}{v\sqrt
n} \biggr\}\stackrel{d}{\longrightarrow} \mathrm{e}^{\mathcal N} .
\]
To finish the proof by Slutsky's theorem it suffices to show
the second factor on the right-hand side of \eqref{eq:rand|M|}
converges to 1
in distribution.
To this end note that
%
\begin{equation}\label{eq:invperp}
S_n:=\sum_{k=1}^nQ_k\prod_{j=1}^{k}\frac1{M_j}=
\sum_{k=1}^n\frac{Q_k}{M_k} \prod_{j=1}^{k-1}\frac{1}{M_j}
\end{equation}
is a perpetuity generated by $(Q/M,1/M)$. Since we are working under the
condition $E{\ln}|M|>0$, we have $E{\ln}|1/M|=-E{\ln}|M|<0$. Furthermore,
by our assumption on $Q$, $E{\ln^+}|Q/M|<\infty$ and thus $(S_n)$
converges in distribution to, say $S$
(see \cite{vervaat}, Theorem~1.6(b)).
Moreover, $P(1/M=0)=0$ and hence by Theorem 1.3 of \cite{air} (see
also \cite{grin} and \cite{blm}, Lemma~2.1) it
follows that the distribution of $S$ is
continuous. In particular, $|S|^{1/v}$ does not have an atom at zero,
which is all that is important for our purposes. Denote by $\nu_n$ and
$\nu$ the distributions
of $|S_n|^{1/v}$ and $|S|^{1/v}$, respectively. We want to show that
$|S_n|^{1/{(v\sqrt{n})}}\stackrel{d}{\longrightarrow}1$.
Consider first an
arbitrary $x\in(0,1)$, fix an arbitrary $m$ and take any $n>m$.
Then
\[
P \bigl(|S_n|^{1/{(v\sqrt{n})}}\le
x \bigr)=\nu_n([0,x^{\sqrt{n}}])\le
\nu_n([0,x^{\sqrt{m}}])\stackrel{n\to\infty}{\longrightarrow}\nu
([0,x^{\sqrt{m}}]).
\]
Letting now $m\to\infty$ we conclude that
$P (|S_n|^{1/{(v\sqrt{n})}}\le x )\to0$ for
$x\in(0,1)$.

Now we take an arbitrary $x\ge1$. Then again we fix some $m$. For
$n>m$ we obtain
\[
P \bigl(|S_n|^{1/{(v\sqrt{n})}}\le x \bigr)=\nu_n([0,x^{\sqrt{n}}])\ge
\nu_n([0,x^{\sqrt{m}}])\stackrel{n\to\infty}{\longrightarrow}\nu
([0,x^{\sqrt{m}}])
\stackrel{m\to\infty}{\longrightarrow}
1
\]
so that $P (|S_n|^{1/{(v\sqrt{n})}}\le x )\to1$ for
$x\ge1$ and it follows that $|S_n|^{1/{(v\sqrt{n})}}\stackrel
d\to1$,
which completes the proof of part (i).

To prove part (ii)
let $\eps_j=\operatorname{sgn}(M_j)$ and consider again
\[
R_n= \Biggl(\prod_{j=1}^nM_j \Biggr)S_n,
\]
where $(S_n)$ is defined by \eqref{eq:invperp}.
By definition of $x^t$,
%
\begin{equation}\label{eq:normRn}
Z_n:=\frac{R_n^{1/{(v\sqrt{n})}}}{\exp(\mu\sqrt n/{v})}=
\operatorname{sgn}(S_n)|S_n|^{1/{(v\sqrt n)}}
\Biggl(\prod_{j=1}^n\eps_j \Biggr)\exp\Biggl(\sum_{j=1}^n\frac{{\ln}|M_j|-\mu
}{v\sqrt{n}} \Biggr).
\end{equation}
For any $m<n$ write the sum in the exponent of \eqref{eq:normRn} as
\[
\sum_{j=1}^m\frac{{\ln}|M_j|-\mu}{v\sqrt{m}}+\sum_{j=m+1}^n\frac
{{\ln}|M_j|-\mu}{v\sqrt n}+ \Biggl(\sqrt{\frac mn}-1 \Biggr)
\sum_{j=1}^m\frac{{\ln}|M_j|-\mu}{v\sqrt{m}}.
\]
Splitting the product of signs on the right-hand side of \eqref
{eq:normRn} in two factors, and using the above equation, we see that
\eqref{eq:normRn} can be written as
%
\begin{equation}\label {eq:Rnm}
Z_n=
\Biggl(\prod_{j=m+1}^n\eps_j \Biggr)Z_m
V_{n,m},
\end{equation}
where
\[
V_{n,m}:=\frac{\operatorname{sgn}(S_n)}{\operatorname
{sgn}(S_m)}\frac{|S_n|^{1/{(v\sqrt
n)}}}{|S_m|^{1/{(v\sqrt m)}}}\exp\Biggl(\sum_{j=m+1}^n\frac{\ln
|M_j|-\mu}{v\sqrt n}+ \Biggl(\sqrt{\frac mn}-1 \Biggr)\sum_{j=1}^m\frac{\ln
|M_j|-\mu}{v\sqrt m} \Biggr).
\]
We claim that $(V_{n,m})$ converges in probability to 1 as long as
$n,m\to\infty$ in such a way that $n-m=\mathrm{o}(n)$. To this end, consider
the first sum in the exponent above. Its variance is $(n-m)/n$, and
thus as long as $n-m=\mathrm{o}(n)$ it goes to 0 in probability by Chebyshev's
inequality.
Furthermore, under the same condition on $n-m$, $\sqrt{\frac
mn}-1=\sqrt{1-\frac{\mathrm{o}(n)}n}-1=\mathrm{o}(1)$\vspace*{-2pt} and thus the second term in the
exponent above goes to 0 in probability as well. (Note that the sum of
$\frac{{\ln}|M_j|-\mu}{v\sqrt m}$ converges in distribution to
$\mathcal
N$ by the classical CLT.)

As for the other factors in $V_{n,m}$, just as in part (i),
$|S_k|^{1/{(v\sqrt k)}} \stackrel P\longrightarrow1$ as \mbox{$k\to\infty$},
where \lq\lq$\stackrel P\longrightarrow$\rq\rq\  denotes
the convergence in probability. It remains to show that
$\operatorname{sgn}(S_n)/\break\operatorname{sgn}(S_m) \stackrel
P\longrightarrow1$. To this end,
note that the equation \eqref{eq:invperp} defining $(S_n)$ holds a.s.
and not only in distribution and thus the sequence $(S_n)$ is a
sequence of partial sums. Hence, by the basic convergence result for
perpetuities (see, e.g., \cite{vervaat}, Lemma~1.7), the conditions
$E{\ln}|1/M|<0$ and $E{\ln^+}|Q/M|<\infty$ suffice for the a.s.
convergence of
the series
\[
\sum_{k=1}^\infty\frac{Q_k}{M_k}\prod_{j=1}^{k-1}\frac1{M_j}.
\]
This, together with the fact that the limit of $(S_m)$ has a continuous
distribution function (see an argument following \eqref{eq:invperp}),
clearly implies that $\operatorname{sgn}(S_n)/\operatorname{sgn}(S_m)
\stackrel P\longrightarrow
1$ (and a.s.) as $m,n\to\infty$.
This proves our claim about $(V_{n,m})$.

We now go back to \eqref{eq:Rnm} and note that the first two random
variables on the right-hand side are independent. Furthermore, the
product $\prod_{j=m+1}^n\eps_j$ converges in distribution to a
symmetric Bernoulli random variable $r$ (see \eqref{eq:symB} above) as
long as $n-m\to\infty$. Hence, if the sequence
\eqref{eq:normRn}
converges in distribution, then its limit, say $Z$, has to satisfy the
distributional identity
$Z\stackrel d= rZ$,
with $r, Z$ independent on the right-hand side. This would complete the
proof since any such $Z$ is symmetric and we would have
%
\begin{equation}\label {eq:Z}Z\stackrel d= r|Z|\stackrel d=
r\mathrm{e}^{\mathcal N},
\end{equation}
where the second equality follows by part (i).

Thus, it remains only to prove that the sequence given in \eqref
{eq:normRn} does, in fact, converge in distribution. We note that this
sequence is tight because the sequence of absolute values
\[
|Z_n|=\frac{\prod_{j=1}^n|M_j|}{\exp(\mu\sqrt
n/{v})}|S_n|^{{1}/{(v\sqrt n)}}
\]
converges in distribution to $\mathrm{e}^{\mathcal N}$ by part (i).
We will complete the proof by showing that every converging subsequence
of \eqref{eq:normRn} converges to the same limit, namely $r\mathrm{e}^{\mathcal
N}$. Let $(k_n)$ be any subsequence for which $(Z_{k_n})$ converges in
distribution. Define
\[
\ell_n=\max\bigl\{m\dvtx k_m\le\sqrt{k_n}\bigr\} \quad \mbox{and set}\quad  m_n=k_n-k_{\ell_n}.
\]
By tightness we can choose a subsequence $(m_{n_j})$ of $(m_n)$ for
which $(Z_{m_{n_j}})$ converges in distribution. Since
$m_{n_j}=k_{n_j}-k_{\ell_{n_j}}$ we have
\[
Z_{k_{n_j}}= \Biggl(\prod_{i=m_{n_j}+1}^{k_{n_j}}\eps_i
\Biggr)Z_{m_{n_j}}V_{k_{n_j},{m_{n_j}}}.
\]
It is readily seen from the construction that
$k_{n_j}-m_{n_j}=\mathrm{o}(k_{n_j})$ and that it converges to infinity.
Therefore, $V_{k_{n_j},m_{n_j}} \stackrel P\longrightarrow1$ as $j\to
\infty$. Thus, if $Z$ and $\overline Z$ are the limits in distribution
of $(Z_{k_n})$ and $(Z_{m_{n_j}})$, respectively, then the above
identity implies that
$Z\stackrel d= r\overline Z$, where $r$ and $\overline Z$ are independent.
It follows that $Z$ is symmetric and thus must satisfy \eqref{eq:Z}.
Since $Z$ was a limit along the arbitrary converging subsequence of
$(Z_n)$, the proof is complete.
\end{pf*}

\section{The case $E{\ln}|M|=0$ and $E|M|>1$}\label{(iii)}
To see the difference between the current situation and the preceding
one, note that because $E{\ln}|M|=0$ the perpetuity defined by (\ref
{eq:invperp}) in the course of the proof of Theorem~\ref{prop:rand|M|}
is not guaranteed to converge. Specifically, consider the following
example: let $X$ be a~non-degenerate, integrable, symmetric random
variable and set
$M=\mathrm{e}^X$.
Then,
\[
E{\ln}|M|=EX=0 \quad \mbox{and}\quad
E |M|>\mathrm{e}^{EX}=1,
\]
where the strict inequality follows from the non-degeneracy of $X$.
By the symmetry of~$X$,
$M$ and $1/M$ have \textit{exactly} the same distribution. Hence, if
we take $Q=1$ in \eqref{perp-it} and \eqref{eq:invperp} we find that
$R_n\stackrel d= S_n$ for $n\ge1$. Therefore, $(R_n)$ and $(S_n)$
converge or diverge in distribution simultaneously. If $M$ is
non-negative this difficulty can be handled
by factoring the largest product in \eqref{perp-it} rather than the
last one; however, the limiting distribution will no longer be
lognormal. It seems very likely that the case of general $M$ is
similar, however it remains open.

Before stating our result in the non-negative case, recall that a
function $h\dvtx\mathbf{R}\to\mathbf R$ is regularly varying at infinity
with index $\rho$ if for every $x>0$
\[
\lim_{t\to\infty}\frac{h(xt)}{h(t)}=x^\rho.
\]
Recall also that any such function can be written as
$h(t)=t^{\rho}\ell(t),$ where $\ell(t)$ is slowly varying at
infinity; that is, for every $x>0$
\[
\lim_{t\to\infty}\frac{\ell(xt)}{\ell(t)}=1.
\]
It follows from that definition in particular that for any $\delta>0$,
$t^\delta\ell(t)\to\infty$ as $t\to\infty$ (see~\cite{bgt} for
these and much more on regularly varying functions).
The following holds true:

\begin{thm}\label{thm:elnm=0}
Let $(Q,M)$ be given by $(\mathrm{e}^Y,\mathrm{e}^X),$ where $EX=0$ and $v^2:=\operatorname{var}(X)$ satisfies $0<v^2<\infty$ (so that $E\ln M=0$; $EM>1$, but may
be infinite; and $v^2=\operatorname{var}(\ln M)<\infty$).
\begin{enumerate}[(iii)]
\item[(i)]
If $EY^2<\infty,$ then, as $n\to\infty,$
%
\begin{equation}\label{elnm=0:sum}
R_n^{1/{(v\sqrt{n})}}
\stackrel d\longrightarrow \mathrm{e}^{|{\mathcal N}|},
\end{equation}
where $|{\mathcal N}|$, the absolute value of the standard normal random
variable, has distribution given by
\[
P(|{\mathcal N}|\le x)=2\Phi(x)-1=\sqrt{\frac2\uppi}\int_0^x\mathrm{e}^{-t^2/2}\,\mathrm{d}t.
\]
\item[(ii)]
Let $h(t):=P(Y> t)$ be a tail function of a random variable $Y$ and
define a~sequen\-ce $(\gamma_n)$ by
\[
\gamma_n:=\inf\biggl\{t\dvtx h(t)\le\frac1n \biggr\},\qquad  n\ge1.
\]
If $h(t)$
is regularly varying at infinity with index $\alpha$
for some $-2<\alpha<0$ then, as $n\to\infty$
%
\begin{equation}\label{elnm=0:max}
R_n^{1/{\gamma_n}}
\stackrel d\longrightarrow \mathrm{e}^{V_\alpha},
\end{equation}
where $V_\alpha$ has Fr\'echet distribution $\Phi_\alpha$ given by
\[
P(V_\alpha\le x)=\Phi_\alpha(x)= \cases{
0,&\quad $x\le0$;\cr
\exp(-x^{\alpha}),&\quad $x>0$.
}
\]
\item[(iii)] Assume that $h(t)$ is regularly varying at infinity with
index $\alpha=-2$, that is, that $h(t)=t^{-2}\ell(t)$, where $\ell
(t)$ is slowly varying at infinity.
If $\lim_{t\to\infty}\ell(t)=\infty$, then \eqref{elnm=0:max}
holds, while if $\lim_{t\to\infty}\ell(t)=0,$ then \eqref
{elnm=0:sum} holds.
\end{enumerate}
\end{thm}

\begin{remarks*}
(1) The only case not covered by the above theorem is $\alpha=-2$ and
$\ell(t)\sim\operatorname{const}$. We suspect that in that case, at least when $X$
and $Y$ are independent, we have
\[
\frac{\max_{1\le k\le n}\{Y_k+S_{k-1}\}}{\sqrt n}\stackrel
d\longrightarrow V_{\alpha}+\sigma{\mathcal N},
\]
where $V_\alpha$ is a Fr\'echet random variable with parameter $\alpha
=-2$, $\sigma$ is the variance of $X$ and $V_\alpha$ and ${\mathcal N}$
are independent, but we have not managed to prove it.

(2) Part \textup{(i)} of Theorem \ref{thm:elnm=0} overlaps with Theorem 2.1\textup{(c)}
of \cite{rs}, where the authors worked under the assumption that $X$
is in the domain of attraction of an $\alpha$-stable law, $1<\alpha
\le2$, but assumed additionally that $X$ and $Y$ are independent.
Under a suitable normalization that agrees with ours for $\alpha=2$
and results in a slightly weaker condition on $Y$ than $EY^2<\infty$,
they obtained as a limit law the distribution of a supremum of the L\'
evy $\alpha$-stable motion on $[0,1]$. A similar comment applies to
the second part of \textup{(iii)}; we assume the regular variation of the tail
of $Y$ but not independence of $X$ and $Y$ as was assumed in \cite
{rs}, Theorem~2.1\textup{(c)}. When $\alpha=2$, the assumptions about the
decay of the tail of $Y$ in both papers are identical.
\end{remarks*}

\begin{pf*}{Proof of Theorem~\ref{thm:elnm=0}}
Set
\[
W_n:=\max_{1\le k\le n}Q_k\prod_{j=1}^{k-1}M_j,
\]
and write
%
\begin{equation}\label{eq:elnm=0}
R_n^{1/{(v\sqrt
{n})}}=W_n^{1/{(v\sqrt{n})}}
\Biggl(\sum_{k=1}^n\frac{Q_k\prod_{j=1}^{k-1}M_j}{W_n} \Biggr)^{1/{(v\sqrt{n})}}.
\end{equation}
Since $x\to\ln x$ is increasing,
\[
W_n^{1/{(v\sqrt{n})}}=\exp\Biggl\{\!\frac1{v\sqrt n}\ln\Biggl(\!\max_{1\le
k\le n}Q_k\prod_{j=1}^{k-1}M_j \!\Biggr)\! \Biggr\}
=\exp\biggl\{\!\frac{\max_{1\le k\le n} \{\ln Q_k+\sum_{j=1}^{k-1}\ln M_j \}
}{v\sqrt n} \!\biggr\}.
\]
We let $Y_k=\ln Q_k$, $X_k=\ln M_k$, $S_k=\sum_{j=1}^kX_k$
and, for any sequence of random variables $(Z_k),$ we will write
$Z_m^*=\max_{1\le k\le m}Z_k$.

Subadditivity of the maxima implies that for any numerical sequences
$(u_k)$, $(w_k)$,
%
\begin{equation}\label{eq:subad}
\max_{1\le k\le n}\{u_k\}-\max
_{1\le k\le n}\{
|w_k|\}\le\max_{1\le k\le n}\{u_k+w_k\}\le\max_{1\le k\le n}\{u_k\}
+\max_{1\le k\le n}\{w_k\}.
\end{equation}
To prove
part (i)
we note that
\[
\frac{|Y_n|^*}{\sqrt n}\stackrel P\longrightarrow0.
\]
Indeed, for $\eps>0$ we have
\[
P \biggl(\frac{|Y_n|^*}{\sqrt n}>\eps\biggr)=1-\bigl(1-P\bigl(|Y|>\eps\sqrt n\bigr)\bigr)^n\le
nP(Y^2>\eps^2n)\to0,\qquad  n\to\infty,
\]
where the last assertion follows from $EY^2<\infty$. Using this and
\eqref{eq:subad} with $u_k=S_{k-1}$ and $w_k=Y_k$ we obtain that
\[
\frac{\max_{1\le k\le n} \{Y_k+S_{k-1} \}}{v\sqrt n}=\frac
{S_{n-1}^*}{v\sqrt n}+\mathrm{o}_P(1),
\]
where $\mathrm{o}_P(1)$ denotes a quantity that goes to zero in probability. Furthermore,
our assumptions on $M_j$'s imply that $(\sum_{j=1}^kX_j)$ is a
sequence of partial sums of random walk whose increments $X_j$ have
mean zero and a finite variance $v^2$. Thus, the Erd\H os--Kac theorem
for the maxima of random walks (see, e.g., \cite{gut}, Theorem~12.2)
implies that
\[
W_n^{1/{(v\sqrt{nv})}}\stackrel d\longrightarrow \mathrm{e}^{|{\mathcal N}|}.
\]
Just as in the proof of Theorem~\ref{prop:rand|M|}, to complete the
proof we need to show that the second factor in (\ref{eq:elnm=0})
converges to 1 in distribution.
But that is clear since, on the one hand,
\[
\Biggl(\sum_{k=1}^n\frac{Q_k\prod_{j=1}^{k-1}M_j}{W_n} \Biggr)^{{1}/{(v\sqrt
{n})}}\le n^{{1}/{(v\sqrt{n})}}=\exp\biggl(\frac{\ln n}{v\sqrt n} \biggr)\to1
\]
and, on the other hand, we clearly have
\[
\sum_{k=1}^n\frac{Q_k\prod_{j=1}^{k-1}M_j}{W_n}\ge1.
\]
This implies that
\[
\sum_{k=1}^n\frac{Q_k\prod_{j=1}^{k-1}M_j}{W_n}\stackrel
P\longrightarrow1,
\]
and proves part (i).

The argument for the second part is parallel to the one just given with
the following adjustment: In the first part, the assumption
$EY^2<\infty$ ensures that the maximum of the random walk with
increments $X_j$, $j<n$, dominates the maximum of $\{Y_k,\  k\le n\}$.
If this assumption is weakened this may no longer be true, and, in
fact, the maximum of $Y_k$'s may dominate. In that case, we can just
use the basics of extreme value theory (Chapter~1 of \cite{llr} being
more than enough) instead of the Erd\H os--Kac theorem to complete the argument.
This time, using \eqref{eq:subad} we write
\[
Y_n^{*}-|S_{n-1}|^*\le\max_{1\le k\le n} \Biggl\{Y_k+\sum_{j=1}^{k-1}X_j \Biggr\}
\le Y_n^*+S_{n-1}^{*}.
\]
By the characterization theorem in the extreme value theory (see, e.g.,
\cite{llr}, Theorem~1.6.2), our assumption on $Y$ is a necessary and
sufficient condition for the existence of constants $(a_n)$, $(b_n)$
for which
\[
a_n(Y^*_n-b_n)\stackrel d\longrightarrow V_\alpha,
\]
where $V_\alpha$ has the Fr\'echet distribution (also referred to as a
type-II extreme value distribution) described above. Furthermore (see
\cite{llr}, Corollary~1.6.3), we may take $a_n=1/\gamma_n$ and
$b_n=0$, $n\ge1$, and, if we do, we obtain that
\[
\frac{Y_n^*}{\gamma_n}\stackrel d\longrightarrow V_\alpha,\qquad  n\to
\infty.
\]
We now observe that $\gamma_n/\sqrt n\to\infty$. In fact, there
exists $\beta>1/2$ such that $\gamma_n\ge n^\beta$ for all
sufficiently large $n$. Indeed, since $h$ is decreasing it is enough to
see that $h(n^\beta)>1/n$. But as $h$ is regularly varying, we have
\[
h(n^\beta)=n^{\beta\alpha}\ell(n^\beta)=\frac{n^{1+\beta\alpha
}\ell(n^\beta)}n=\frac{(n^\beta)^{1/\beta+\alpha}\ell
(n^\beta)}n=\frac{(n^\beta)^{\delta}\ell(n^\beta)}n.
\]
If we now take $1/2<\beta<-1/\alpha$ (which is possible since
$-2<\alpha<0$) then $\delta>0$ and the numerator on the right-hand
side goes to infinity with $n$, proving the claim that $h(n^\beta
)>1/n$ for large $n$. We now get
\[
\frac{|S_{n-1}|^*}{\gamma_n}=\frac{|S_{n-1}|^*}{\sqrt n}\frac{\sqrt
n}{\gamma_n}\stackrel P\longrightarrow0,
\]
and also, since $\gamma_n/\sqrt n\to\infty$,
\[
1\le\Biggl(\sum_{k=1}^n\frac{Q_k\prod_{j=1}^{k-1}M_j}{W_n} \Biggr)^{
{1}/{\gamma_n}}\le n^{{1}/{\gamma_n}}\to1,
\]
which proves the second part.

Finally, the last part follows by essentially the same reasoning.
Assume $\ell(t)\to\infty$ as $t\to\infty$. By the just-given
argument, to establish \eqref{elnm=0:max}, it suffices to verify
$\gamma_n/\sqrt n\to\infty$ as $n\to\infty$. Assume this is not
the case. Then there exist $C<\infty$ and an infinite subsequence
$(n_k)$ such that $\gamma_{n_k}/\sqrt{n_k}\le C$ for all $k\ge1$. By
the definition of $(\gamma_n)$ this means that $P(Y\ge C\sqrt
{n_k})\le1/n_k$, that is, that $\ell(C\sqrt{n_k})\le C^2$. But that
contradicts the assumption that $\ell(t)\to\infty$ as $t\to\infty
$. If, on the other hand, $\ell(t)\to0$ as $t\to\infty$, then
$\gamma_n/\sqrt n\to0$ (for otherwise there would exist $c>0$ and a
subsequence $(n_k)$ such that $\gamma_{n_k}/\sqrt{n_k}\ge c$, $k\ge
1$, implying that $\ell(c\sqrt{n_k})\ge c^2$ and contradicting $\ell
(t)\to0$ as $t\to\infty$). Now, it follows from the proof of the
first part that $\gamma_n/\sqrt n\to0$ is enough to conclude \eqref
{elnm=0:sum}. The proof of part (iii) is completed.
\end{pf*}

\section{The case $E{\ln}|M|=0$ and $E|M|=1$}
Under this assumption we have that $|M|\equiv1$, that is, $M$ takes on values
$\pm1$. Consequently, because of a non-degeneracy assumption on $M$,
its expected value must satisfy $-1<EM<1$. We have
\begin{thm}\label {prop:clt} Suppose that $EQ^2<\infty$. Then, as
$n\to
\infty$,
%
\begin{equation}\label {eq:clt}\frac{R_n}{\sqrt{n}}\stackrel
{d}{\longrightarrow}
\beta{\mathcal N},
\end{equation}
where $\beta^2=EQ^2+2\frac{EQ}{1-EM}E(QM)$ and $\mathcal N$ is the
standard normal random variable.
\end{thm}

\begin{remark*}
Since $-1<EM<1$, by straightforward calculation
we see that $ER_n=\mathrm{O}(1)$ and
$\operatorname{var}(R_n)=\beta^2n+\mathrm{O}(1)$. Thus, (\ref{eq:clt}) is
equivalent to
\[
\frac{R_n-ER_n}{\sqrt{\operatorname{var}(R_n)}}\stackrel
d\longrightarrow{\mathcal N}.
\]
\end{remark*}

\begin{pf}
Let $\alpha:=EM$ and $q:=EQ$ so that
\[
\beta^2=EQ^2+2\frac q{1-\alpha}EQM=\operatorname{var}\biggl(Q+\frac q{1-\alpha}M\biggr).
\]
To prove (\ref{eq:clt}) we write
%
\begin{equation}\label {eq:mart1}R_n=\sum_{k=1}^nQ_k\prod
_{j=1}^{k-1}M_j=\sum
_{k=1}^n(Q_k-q)\prod_{j=1}^{k-1}M_j+q\sum_{k=1}^n\prod_{j=1}^{k-1}M_j.
\end{equation}
Furthermore, as
\[
\sum_{k=1}^n(M_k-\alpha)\prod_{j=1}^{k-1}M_j=-\alpha+\prod
_{j=1}^nM_j+(1-\alpha)\sum_{k=2}^n\prod_{j=1}^{k-1}M_j
=-1+\prod_{j=1}^nM_j+(1-\alpha)\sum_{k=1}^n\prod_{j=1}^{k-1}M_j,
\]
the second term on the right-hand side of (\ref{eq:mart1}) is
\[
\frac q{1-\alpha}\sum_{k=1}^n(M_k-\alpha)\prod_{j=1}^{k-1}M_j-\frac
q{1-\alpha}\prod_{j=1}^nM_j+\frac q{1-\alpha}.
\]
Set
\[
d_k:=\biggl(Q_k-q+\frac q{1-\alpha}(M_k-\alpha)\biggr)\prod_{j=1}^{k-1}M_j,\qquad
k=1,\dots, n.
\]
Then, by independence of $(Q_j,M_j)$, $j\ge1$, $(d_k)$ is a martingale
difference sequence with respect to $({\mathcal F}_n)$ where
${\mathcal
F}_k=\sigma(M_1,Q_1,\dots, M_k,Q_k)$, $k\ge1$. Let $E_k$ denote the
conditional expectation given ${\mathcal F}_k$. Since $M_j^2=1$ we have
\[
E_{k-1}d_k^2=
E \biggl(Q-q+\frac q{1-\alpha}(M-\alpha) \biggr)^2=\beta^2,
\]
so that, trivially,
\[
\frac{\sum_{k=1}^nE_{k-1}d_k^2}{n}\stackrel P\longrightarrow\beta^2.
\]
Moreover, since the $M_k$'s are uniformly
bounded, for a given $\eps>0$ and $n$ sufficiently large
\[
E_{k-1}d_k^2I_{\{|d_k|>\eps\sqrt n\}}\le4EQ^2I_{\{|Q|>\eps\sqrt n/2\}}.
\]
Since $EQ^2<\infty$, the last quantity converges to zero as $n\to
\infty$ by the dominated convergence theorem. This verifies the
conditional version of Lindeberg's condition:
\[
\forall\eps>0\qquad  \frac{\sum_{k=1}^nE_{k-1}d_k^2I_{\{|d_k|>\eps\sqrt
n\}}}{n}\stackrel P\longrightarrow0.
\]
It follows by the martingale version of the CLT (see, e.g.,
\cite{bil}, Theorem~35.12) that
\[
\frac{\sum_{k=1}^nd_k}{\sqrt{n}}\stackrel d\longrightarrow\beta
{\mathcal N}.
\]
Now, by the above manipulations we have
\[
\frac{R_n}{\sqrt n}=
\frac{\sum_{k=1}^nd_k}{\sqrt n}+\frac{q}{(1-\alpha)\sqrt
n}-\frac q{(1-\alpha)\sqrt n}\prod_{j=1}^{n}M_j.
\]
Since each of
the last two terms goes to 0 (deterministically and in
probability, respectively), Theorem~\ref{prop:clt} follows.
\end{pf}

\section*{Acknowledgements}
The first author is supported in part by the NSA Grant
\#H98230-09-1-0062 and would like to thank S.~Kwapie\'n for several
useful conversations on the subject of this paper. We would also like
to thank A.~Iksanov and D.~Buraczewski for pointing us to some work on
a divergent case, including \cite{rs} and \cite{pakes}.

\printhistory

\end{document}